\documentclass[12pt]{article}
\usepackage{color}
\definecolor{darkblue}{rgb}{0.00,0.05,0.50}
\usepackage[colorlinks,filecolor=blue,citecolor=darkblue]{hyperref}

\usepackage{xcolor}
\usepackage{hyperref}

\definecolor{linkcolor}{HTML}{120a8F} % колір посилань
\definecolor{urlcolor}{HTML}{120a8F} % колір гіперпосилань

\hypersetup{pdfstartview=FitH,  linkcolor=linkcolor,urlcolor=urlcolor, colorlinks=true}

\usepackage{amsfonts,amssymb,amsmath,amsthm}
\usepackage{mathrsfs}
\usepackage{url}
\usepackage{enumerate}
\usepackage[ukrainian, russian, english]{babel}
\usepackage[cp1251]{inputenc}
\usepackage{amscd, latexsym,cite,verbatim,texdraw,floatflt, caption2,pb-diagram} %new

\theoremstyle{plain}
\newtheorem{theorem}{Теорема}[section]

\newtheorem{proposition}{Твердження}[section]
\newtheorem{corollary}{Наслідок}[section]

\theoremstyle{definition}

\newtheorem{remark}{Зауваження}[section]
\newcommand{\keywords}{\textbf{Key words and phrases: }\medskip}
\newcommand{\subjclass}{\textbf{MSC 2020: }\medskip}
\renewcommand{\abstract}{\textbf{Abstract }\medskip}

\numberwithin{equation}{section}

\begin{document} \selectlanguage{ukrainian}
\thispagestyle{empty}

\title{}

UDC 517.51 \vskip 5mm

\begin{center}
\textbf{\Large Approximation of functions  of many variables from the generalization Nikol'skіі--Besov type classes in the uniform and integral metrics}
\end{center}

\vskip 3mm

\begin{center}
\textbf{\Large  Наближення функцій багатьох змінних із узагальнених класів
типу Нікольського--Бесова у рівномірній та інтегральній метриках}
% цілими функціями експоненціального типу
\end{center}
\vskip0.5cm

\begin{center}
\large{M.\,I.~Gromyak,  O.\,Ya.~Radchenko, S.\,Ya.~Yanchenko}
\end{center}
\begin{center}
\large{М.\,І.~Гром'як, О.\,Я.~Радченко, C.\,Я.~Янченко}
\end{center}
\vskip0.5cm

\begin{abstract}

We obtain the exact order estimates of the approximation of the functions of many variables from the generalized  Nikol'skіі--Besov  classes $B^{\Omega}_{p,\theta}(\mathbb{R}^d)$ by sums of de la Vallee Poussin type in the metrics space $L_{\infty}(\mathbb{R}^d)$ and $L_{1}(\mathbb{R}^d)$. These classes of functions for some given $\Omega$ coincide with the well-known classical  Nikol'skіі--Besov isotropic classes.

\vskip 3 mm

Одержано точні за порядком оцінки наближення функцій багатьох змінних із узагальнених класів Нікольського--Бєсова $B^{\Omega}_{p,\theta}(\mathbb{R}^d)$, які при деякому $\Omega$ збігаються з відомими класичними ізотропними класами Нікольського--Бєсова,  сумами Валле Пуссена в метриках просторів $L_{\infty}(\mathbb{R}^d)$ і $L_{1}(\mathbb{R}^d)$.
\end{abstract}

\vskip 5 mm

\subjclass{41A30, 41A50, 41A63, 42A38.}

\keywords{Nikol'skii--Besov type classes, entire functions of exponen\-tial type, Fourier transform, de la Vallee Poussin sums.}

\vskip 0.7 cm

%%%%%% Вступ %%%%%%%%%%%%%%%%%%%%%%%%%%%%%%%%%%%%%%%%%%%%%%%%%%%%%%%%%%%%%%%%%%%%%%%%

\section{Вступ}

У статті  продовжується досліджуватися питання наближення класів $B^{\Omega}_{p,\theta}(\mathbb{R}^d)$~\cite{Muronyuk_2014_UMG, Muronyuk_2013_zb} функцій багатьох змінних у просторах $L_{\infty}(\mathbb{R}^d)$ та $L_{1}(\mathbb{R}^d)$ за допомогою цілих функцій експоненціального типу (див., наприклад, \cite[гл.~3]{Nikolsky_1969_book}). Дані класи функцій при  $\Omega(t)=t^r$ збігаються з добре відомими ізотропними класами Нікольського--Бєсова  $B^{r}_{p,\theta}(\mathbb{R}^d)$~\cite{Besov_1961, Nikolsky_1951}. Деякі з отриманих результати були анонсовані в \cite{Yanchenko-UMJ-2016}, тут ми наводимо їх з повним доведенням та усіма необхідними коментарями.

\section{Основні позначення, означення класів та апроксимаційних характеристик}

Нехай $\mathbb{R}^d$, $d\geqslant 1$~--- $d$-вимірний евклідів простір з елементами $\boldsymbol{x}=(x_1,\ldots,x_d)$, ${(\boldsymbol{x},\boldsymbol{y})=x_1y_1+\ldots+x_dy_d}$. Через ${L_p(\mathbb{R}^d)}$ позначимо простір
вимірних на $\mathbb{R}^d$ і сумовних у степені $p$, $1\leqslant p < \infty$, відповідно суттєво обмежених при $p=\infty$, функцій $f(\boldsymbol{x})=f(x_1,\ldots,x_d)$ зі скінченною нормою, яка визначається таким чином:
$$
\|f\|_{L_p(\mathbb{R}^d)}:=
\left(\int\limits_{\mathbb{R}^{d}}|f(\boldsymbol{x})|^{p}d\boldsymbol{x}
\right) ^{\frac{1}{p}}, \ 1\leqslant p<\infty,
$$
$$
\|f\|_{L_{\infty}(\mathbb{R}^d)}:=\mathop {\rm ess \sup}\limits_{\boldsymbol{x}\in \mathbb{R}^d}
  |f(\boldsymbol{x})|.
 $$

Для функції $f\in L_p(\mathbb{R}^d)$ у точці $\boldsymbol{x}$ з кроком $\boldsymbol{h}\in \mathbb{R}^d$ позначимо через $\Delta_{\boldsymbol{h}}^l f(\boldsymbol{x})$ кратну різницю порядку $l \in \mathbb{N}$
$$
 \Delta_{\boldsymbol{h}}^l f(\boldsymbol{x})=
 \Delta_{\boldsymbol{h}}\Delta_{\boldsymbol{h}}^{l-1}f(\boldsymbol{x}), \ \  \Delta_{\boldsymbol{h}}^0 f(\boldsymbol{x})=f(\boldsymbol{x})
$$
де $\Delta_{\boldsymbol{h}} f(\boldsymbol{x})=f(\boldsymbol{x}+\boldsymbol{h})-f(\boldsymbol{x})$.

Кратну різницю $\Delta_{\boldsymbol{h}}^l f(\boldsymbol{x})$ також можна записати у вигляді
$$
\Delta_{\boldsymbol{h}}^l f(\boldsymbol{x})=\sum\limits_{j=0}^l (-1)^{j+l}C^j_l f(\boldsymbol{x}+l\boldsymbol{h}),
$$
де $C^{j}_l$~--- біноміальні коефіцієнти.

На основі кратної різниці $\Delta_{\boldsymbol{h}}^l f(\boldsymbol{x})$ означимо модуль неперервності $l$-го порядку функції $f\in L_p(\mathbb{R}^d)$, який будемо позначати $\Omega_l(f,t)_p$, такою формулою:
$$
\Omega_l(f,t)_p= \sup \limits_{|\boldsymbol{h}|\leqslant t} \|\Delta_{\boldsymbol{h}}^l f(\cdot)\|_{L_p(\mathbb{R}^d)},
$$
де $|\boldsymbol{h}|=\sqrt{h_1^2+\ldots+h_d^2}$~--- евклідова норма вектора $\boldsymbol{h}$.

Нехай $\Omega(t):=\Omega$~--- функція типу модуля неперервності порядку $l$, тобто функція, яка визначена і неперервна на $\mathbb{R_+}$, що задовольняє такі умови:

1) $\Omega(0)=0$, $\Omega(t)> 0$ для $t> 0$;

2) $\Omega(t)$ неспадна на $\mathbb{R_+}$;

3) для всіх $n\in \mathbb{Z_+}$,  $ \Omega(nt)\leqslant C_1 n^l \Omega(t)$, де $C_1>0$ не залежить від $n$ i $t$.

Множину таких функцій $\Omega$ позначимо через $\Psi_{l}$.

Додатково будемо вимагати, щоб функція $\Omega$ задовольняла умови $(S^{\alpha})$ та $(S_l)$. У  літературі ці умови  називають умовами Барi--Стєчкiна~\cite[\S2]{Bari_Stechkin} і вони описуються в термінах двох понять (майже зростання та майже спадання функцій), які запроваджені С.\,Н.~Бернштейном~\cite[\S97, \S108]{Bernstein}. Сформулюємо ці умови:
\begin{itemize}
\item[І)] функція однієї змінної $\varphi(\tau)\geqslant 0$ задовольняє
умову $(S^{\alpha})$, якщо існує таке $\alpha>0$, що $\varphi(\tau)/\tau^{\alpha}$  майже
зростає, тобто існує така незалежна від $\tau_1$  i
$\tau_2 $ стала ${C_2>0}$, що
 $$
 \frac{\varphi(\tau_1)}{\tau_1^{\alpha}} \leqslant C_2\frac{\varphi(\tau_2)}{\tau_2^{\alpha}},
  \ \ \ \ 0<\tau_1 \leqslant \tau_2 \leqslant 1;
 $$

\item[ІІ)] функція однієї змінної $\varphi(\tau)\geqslant 0$ задовольняє
умову $(S_l)$, якщо  існує таке $\gamma$, $0<\gamma<l$, що
$\varphi(\tau)/\tau^{l-\gamma}$ майже спадає, тобто
існує така незалежна від $\tau_1$ i $\tau_2 $ стала $C_3>0$,  що
$$
\frac{\varphi(\tau_1)}{\tau_1^{l-\gamma}} \geqslant
C_3\frac{\varphi(\tau_2)}{\tau_2^{l-\gamma}},  \ \ \ \ 0<\tau_1
\leqslant \tau_2 \leqslant 1 .
$$
\end{itemize}

У випадку, коли для $\Omega$ виконано умову $(S^{\alpha})$, будемо говорити, що $\Omega$ належить множині $S^{\alpha}$, а якщо умову $(S_l)$, то~--- множині $S_l$. Якщо $\Omega$ одночасно належить множинам $\Psi_{l}$, $S^{\alpha}$ і $(S_l)$, використовуватимемо запис $\Omega \in \Phi_{\alpha, l}$, тобто множина $\Phi_{\alpha, l}$ визначається співвідношенням $\Phi_{\alpha, l}=\Psi_l\cap S^{\alpha} \cap S_l$.

Зазначимо, що до множини $\Phi_{\alpha, l}$ належать, наприклад, функції
$$
\Omega(t)=
\begin{cases}
    t^r \left( \log_{\,2} ^+ \frac{1}{t} \right)^\beta, & t>0, \\
    0, & t=0,
 \end{cases}
$$
де $\log_{\,2}^+ t= \max\{1, \log_{\,2} t\}$, $\alpha<r<l$, a $\beta$~--
фіксоване дійсне число.

Для $1 \leqslant p, \theta \leqslant\infty$ і заданої функції $\Omega \in \Psi_l$ простори $B^{\Omega}_{p,\theta}(\mathbb{R}^d)$ означаються таким чином \cite{Muronyuk_2014_UMG, Muronyuk_2013_zb}:
$$
B^{\Omega}_{p,\theta}(\mathbb{R}^d):=\left\{f \in L_p(\mathbb{R}^d) \colon
\|f\|_{B^{\Omega}_{p,\theta}(\mathbb{R}^d)} < \infty \right\},
$$
де норма задається співвідношеннями
\begin{equation}\label{B-Omega-norm1}
\|f\|_{B^{\Omega}_{p,\theta}(\mathbb{R}^d)} = \|f\|_{L_p(\mathbb{R}^d)} +
 \left(\int\limits_0^{\infty} \left(\frac{\Omega_l(f,t)_p}{\Omega(t)}\right)^{\theta}\frac{dt}{t}\right)^{\frac{1}{\theta}}, \ \ 1 \leqslant \theta < \infty,
\end{equation}
і
\begin{equation}\label{B-Omega-norm2}
\|f\|_{B^{\Omega}_{p,\infty}(\mathbb{R}^d)} = \|f\|_{L_p(\mathbb{R}^d)} + \sup\limits_{t>0}\ \frac{\Omega_l(f,t)_p}{\Omega(t)}
\end{equation}
при $\theta=\infty$.

Простори $B^{\Omega}_{p,\theta}(\mathbb{R}^d)$ з так введеною нормою будуть банаховими.

Зауважимо, що простори $B^{\Omega}_{p,\theta}(\mathbb{R}^d)$, як уже зазначалося, у випадку $\Omega(t)=t^r$, $0<r<l$, співпадають з просторами $B^{r}_{p,\theta}(\mathbb{R}^d)$. Нагадаємо також, що простори $H^r_p(\mathbb{R}^d)\equiv B^r_{p,\infty}(\mathbb{R}^d)$ уперше були розглянуті С.\,М.~Нікольським~\cite{Nikolsky_1951}, а простори $B^r_{p,\theta}(\mathbb{R}^d)$ були введені О.\,В.~Бєсовим~\cite{Besov_1961}. Таким чином, простори
$B^{\Omega}_{p,\theta}(\mathbb{R}^d)$ є узагальненням за гладкісним параметром відомих ізотропних просторів Нікольського--Бєсова.

Також зазначимо, що у випадку $d=1$ простори $B^{\Omega}_{p,\theta}(\mathbb{R})$ збігаються з просторами $S^{\Omega}_{p,\theta}B(\mathbb{R})$ (див., наприклад, \cite{Stasuk-Yanchenko-Anal-math, Yanchenko-UMJ-2016}), які є узагальненням відомих просторів функцій з домінуючою мішаною похідною. У свою чергу, при певному заданні параметра $\Omega$ простори $S^{\Omega}_{p,\theta}B(\mathbb{R}^d)$, $d\geqslant 1$, збігаються з просторами  $S^{\boldsymbol{r}}_{p,\theta}B(\mathbb{R}^d)$, які  були введені С.\,М.~Нікольським~\cite{Nikolsky_63} при $\theta=\infty$ та Т.\,І.~Амановим~\cite{Amanov_1965} при $1\leqslant \theta < \infty$ і з точки зору знаходження точних за порядком оцінок деяких апроксимаційних характеристик досліджувалися, зокрема,  в~\cite{WangHeping-SunYongsheng-1995, WangHeping-SunYongsheng-1999-AppT, Yanchenko_CMP2020, Yanchenko_YMG_2019}.

Наведемо один результат, що встановлений в \cite{Muronyuk_2014_UMG}, який дає можливість означити норму функцій з просторів $B^{\Omega}_{p,\theta}(\mathbb{R}^d)$ в іншій формі. В подальшому це зумовлює використання перетворення Фур'є у теорії даних просторів із застосуванням узагальнених функцій  (див., наприклад, \cite[гл.~11]{Berezanskyi2014}) \cite[гл.~2]{Vladimirov}, \cite{Lizorkin_69}, \cite[гл.~1, \S5]{Nikolsky_1969_book}. Для цього додатково наведемо необхідне позначення  і означення.

Нехай $S=S(\mathbb{R}^d)$~--- простір Л.~Шварца основних нескінченно диференційовних на $\mathbb{R}^d$ комплекснозначних функцій $\varphi$, які спадають на нескінченності разом зі своїми похідними
швидше за будь-який степінь функції $|\boldsymbol{x}|^{-1}$, %$\left(x_1^2+\ldots+x_d^2\right)^{-\frac{1}{2}}$
що розглядається з відповідною топологією. Через $S'$ позначимо простір лінійних неперервних функціоналів над $S$. Зазначимо, що елементами простору $S'$ є узагальнені функції. Якщо $f\in S'$, $\varphi\in S$, то $\langle f,\varphi \rangle$  позначає значення $f$ на $\varphi$.

Перетворення Фур'є $\mathfrak{F}\varphi\colon S\rightarrow S$ визначається за формулою
$$
(\mathfrak{F}\varphi)(\boldsymbol{\lambda})=\frac{1}{(2\pi)^{\frac{d}{2}}}\int
\limits_{\mathbb{R}^d}\varphi(\boldsymbol{t}) e^{-i(\boldsymbol{\lambda},\boldsymbol{t})}d\boldsymbol{t}\equiv
\widetilde{\varphi}(\boldsymbol{\lambda}),
$$
де $\boldsymbol{\lambda}=(\lambda_1,\ldots,\lambda_d)\in \mathbb{R}^d$,
$\boldsymbol{t}=(t_1,\ldots,t_d)\in \mathbb{R}^d$ i $(\boldsymbol{\lambda},\boldsymbol{t})=\sum
\limits_{i=1}^d \lambda_it_i$~--- скалярний добуток в $\mathbb{R}^d$
векторів $\boldsymbol{\lambda}$ i $\boldsymbol{t}$.

Обернене перетворення Фур'є задається таким чином:
$$
(\mathfrak{F}^{-1}\varphi)(\boldsymbol{t})=\frac{1}{(2\pi)^{\frac{d}{2}}}\int
\limits_{\mathbb{R}^d}\varphi(\boldsymbol{\lambda})
e^{i(\boldsymbol{\lambda},\boldsymbol{t})}d\boldsymbol{\lambda}\equiv \widehat{\varphi}(\boldsymbol{t}).
$$

Перетворення Фур'є узагальнених функцій (для нього ми зберігаємо те
ж позначення) визначається згідно з формулою
$$
\langle \mathfrak{F}f,\varphi\rangle=\langle f,\mathfrak{F}\varphi
\rangle, \ \ \ \langle \widetilde{f},\varphi\rangle=\langle
f,\widetilde{\varphi} \rangle,
$$
де $f\in S'$, a $\varphi \in S$.

Обернене перетворення узагальнених функцій  також позначимо
$\mathfrak{F}^{-1}f$, і визначається воно аналогічно до прямого
перетворення Фур'є за правилом
$$
\langle \mathfrak{F}^{-1}f,\varphi\rangle=\langle
f,\mathfrak{F}^{-1}\varphi \rangle, \ \ \ \langle
\widehat{f},\varphi\rangle=\langle f,\widehat{\varphi} \rangle.
$$

Носієм неперервної на $\mathbb{R}^d$ функції $\varphi$ називається замикання множини точок ${\boldsymbol{x} \in \mathbb{R}^d}$, де $\varphi(\boldsymbol{x})\neq 0$, і позначається $\mbox{supp}\,\varphi$.

Узагальнена функція  $f$ перетворюється в нуль на відкритій множині $G$, якщо $\langle f,\varphi \rangle=0$ для всіх $\varphi \in S$ і $\mbox{supp}\,\varphi\subset G$. Об'єднання всіх околів, у яких $f$ перетворюється в нуль є відкритою множиною, яку називають нульовою множиною узагальненої функції $f$ і позначають $G_f$. Носієм узагальненої функції називають доповнення множини $G_f$ до $\mathbb{R}^d$, тобто замкнену множину $\mbox{supp}\,f=\bar{G}_f$.

Зазначимо, що кожна функція $f \in L_p(\mathbb{R}^d)$, $1\leqslant p \leqslant \infty$, визначає лінійний неперервний функціонал на $S$ згідно з формулою
$$
\langle f,\varphi \rangle = \int \limits_{\mathbb{R}^d} f(\boldsymbol{x})\varphi(\boldsymbol{x}) d\boldsymbol{x}, \ \ \varphi\in S,
$$
і, як наслідок, у цьому сенсі вона є елементом $S'$. Тому перетворення Фур'є функції ${f \in L_p(\mathbb{R}^d)}$, $1\leqslant p \leqslant \infty$, можна
розглядати, як перетворення Фур'є узагальненої функції $\langle f,\varphi \rangle$.

Далі для $s\in \mathbb{Z}_+$ позначимо
$$
A_{2^s}=\Big\{\boldsymbol{\lambda}=(\lambda_1, \ldots, \lambda_d)\in \mathbb{R}^d \colon -2^{s}< \lambda_j < 2^{s}, \quad j=\overline{1,d}\Big\}
$$
i
$$
G_q(A_{2^n})=\Big\{g\in L_q(\mathbb{R}^d)\colon \mbox{supp}\,
\mathfrak{F}g\subseteq A_{2^n} \Big\}, \quad\quad n\in \mathbb{Z_+},
$$
---~множина цілих функцій експоненціального типу, які належать $L_q(\mathbb{R}^d)$ і носій перетворення Фур'є яких міститься в $A_{2^n}$. Тоді для $f\in L_q(\mathbb{R}^d)$ покладемо
$$
E_{A_{2^n}}\big(f\big)_{L_q(\mathbb{R}^d)}:=E_{n}\big(f\big)_{L_q(\mathbb{R}^d)}:=\inf \limits_{g\in G_q(A_{2^n})}\|f(\cdot)-g(\cdot)\|_{L_q(\mathbb{R}^d)}.
$$
Дана величина називається найкращим наближенням функції $f$ в метриці простору $L_q(\mathbb{R}^d)$ функціями з $G_q(A_{2^n})$. Відповідно для функціонального класу $F\subset L_q(\mathbb{R}^d)$
покладемо
\begin{equation}\label{En-F}
E_n\big(F\big)_{L_q(\mathbb{R}^d)}:=\sup \limits_{f\in F}E_n\big(f\big)_{L_q(\mathbb{R}^d)}.
\end{equation}

%\textbf{ Теорема~А}~\cite{Muronyuk_2014_UMG}.
\begin{proposition}[\!\!\cite{Muronyuk_2014_UMG}] \label{T-decompoz}
Нехай $1\leqslant p, \theta \leqslant \infty$ і $\Omega \in
\Phi_{\alpha, l}$, $\alpha>0$, $l\in \mathbb{N}$. Функція $f$ належить простору $B^{\Omega}_{p,\theta}(\mathbb{R}^d)$, тоді і тільки тоді, коли вона зображується збіжним у метриці $L_p(\mathbb{R}^d)$ рядом
\begin{equation}\label{ef-Qs}
f(\boldsymbol{x})=\sum\limits_{s=0}^{\infty}Q_{\boldsymbol{s}}(\boldsymbol{x}), \ Q_{\boldsymbol{s}}(\boldsymbol{x})= Q_{s,\ldots,s}(\boldsymbol{x})
\end{equation}
де $Q_{\boldsymbol{s}}(\boldsymbol{x})$~--- цілі функції експоненціального типу степенів не вище $2^s$ по кожній змінній, для яких виконуються умови
$$
\left(\sum\limits_{s=0}^{\infty} \big(\Omega(2^{-s})\big)^{-\theta} \|Q_{\boldsymbol{s}}(\cdot)\|^{\theta}_{L_p(\mathbb{R}^d)}\right)^{\frac{1}{\theta}}<\infty,
\  \mbox{якщо} \  1\leqslant \theta<\infty,
$$
і
$$
\sup\limits_{s\in \mathbb{Z}_+}\frac{\|Q_{\boldsymbol{s}}(\cdot)\|_{L_p(\mathbb{R}^d)}}{\Omega(2^{-s})}<\infty, \ \mbox{якщо}  \  \theta=\infty.
$$
Більше того, мають місце співвідношення
$$
\|f\|_{B^{\Omega}_{p,\theta}(\mathbb{R}^d)}\ll \left(\sum\limits_{s=0}^{\infty} \big(\Omega(2^{-s}))^{-\theta} \|Q_{\boldsymbol{s}}(\cdot)\|^{\theta}_{L_p(\mathbb{R}^d)}\right)^{\frac{1}{\theta}}
\ \mbox{при} \ \  1\leqslant \theta<\infty
$$
і
$$
\|f\|_{B^{\Omega}_{p,\theta}(\mathbb{R}^d)}\ll \sup\limits_{s\in \mathbb{Z}_+}\frac{\|Q_{\boldsymbol{s}}(\cdot)\|_{L_p(\mathbb{R}^d)}}{\Omega(2^{-s})} \ \mbox{при} \  \theta=\infty.
$$

Якщо ж, крім цього, частинні суми $n$-го порядку ряду \eqref{ef-Qs} реалізують найкраще наближення~\eqref{En-F} або принаймні дають порядок найкращого наближення, то
\begin{equation}\label{f-norm-dek1}
\|f\|_{B^{\Omega}_{p,\theta}(\mathbb{R}^d)}\asymp \left(\sum\limits_{s=0}^{\infty} \big(\Omega(2^{-s}))^{-\theta} \|Q_{\boldsymbol{s}}(\cdot)\|^{\theta}_{L_p(\mathbb{R}^d)}\right)^{\frac{1}{\theta}}
\ \ \mbox{при} \ \  1\leqslant \theta<\infty
\end{equation}
і
\begin{equation}\label{f-norm-dek2}
\|f\|_{B^{\Omega}_{p,\theta}(\mathbb{R}^d)}\asymp \sup\limits_{s\in \mathbb{Z}_+}\frac{\|Q_{\boldsymbol{s}}(\cdot)\|_{L_p(\mathbb{R}^d)}}{\Omega(2^{-s})} \ \ \mbox{при} \  \theta=\infty.
\end{equation}
\end{proposition}

Одразу зауважимо, що для просторів $B^r_{p,\theta}(\mathbb{R}^d)$ твердження аналогічне твердженню~\ref{T-decompoz} встановив П.\,І.~Лізоркін~\cite{Lizorkin_1968_sib}. Дослідження ізотропних класів Нікольського--Бєсова та їхніх узагальнень також проводилися у роботах \cite{LiuYongping-XuGuiqiao-2002, JiangYanjie-LiuYongping-2000}.

Тут і далі по тексту для додатних величин $A$ і  $B$ використовуємо запис  $A\asymp B$, який означає, що існують такі додатні  сталі $C_4$ і $C_5$, які не залежать від одного істотного параметра у величинах  $A$ і  $B$ (наприклад, у співвідношеннях \eqref{f-norm-dek1} і \eqref{f-norm-dek2}~--- від функції $f$), що ${C_4 A \leqslant B \leqslant C_5 A}$. Якщо тільки $B\leqslant C_5 A $ $\big(B \geqslant C_4 A\big)$, то пишемо
$B\ll A$ ${\big(B \gg A \big)}$. Всі сталі $C_i$, $i=1,2,\dots$, які зустрічаються у роботі, можуть залежати лише від параметрів, що входять в означення класу, метрики, в якій оцінюється похибка наближення, та розмірності простору $\mathbb{R}^d$.

На основі  твердження~\ref{T-decompoz} дамо еквівалентне означення норми функцій із просторів $B^{\Omega}_{p,\theta}(\mathbb{R}^d)$, яким будемо користуватися у подальших міркуваннях.

Розглянемо неперіодичний аналог ядра Валле Пуссена~\cite[\S 8.6]{Nikolsky_1969_book} (див. також \cite{Nikolskii_1968_sib})

$$
 V_{2^s}(\boldsymbol{x})=\frac{1}{2^{sd}}\prod
 \limits_{j=1}^{d}\frac{\cos{2^s x_j}-\cos{2^{s+1}x_j}}{x_j^2}, \ \
 j=\overline{1,d}, \ \ s\in \mathbb{N}\cup\{0\}.
$$

Дане ядро має такі властивості:

\ 1) $V_{2^s}(\boldsymbol{z})=V_{2^s}(z_1,\dots,z_d)$~--- ціла функція
експоненціального типу степеня $2^{s+1}$ за кожною змінною $z_j$, $j=\overline{1,d}$,
обмежена і сумовна на $\mathbb{R}^d$;

\ 2) $\left(\frac{2}{\pi}\right)^{\frac{d}{2}}\tilde{V}_{2^s}=\frac{1}{\pi^d}
\int_{\Box_{2^s}} V_{2^s}(\boldsymbol{t}) e^{-i\boldsymbol{t}\boldsymbol{x}}d\boldsymbol{t}$, де $\Box_{2^s}=\big\{ |x_j|\leqslant 2^s, \ j=\overline{1,d}\big\}$;

\ 3) $\frac{1}{\pi^d}\int_{\mathbb{R}^d}V_{2^s}(\boldsymbol{t})d\boldsymbol{t}=1$;

\ 4) $\frac{1}{\pi^d}\int_{\mathbb{R}^d}|V_{2^s}(\boldsymbol{t})|d\boldsymbol{t}\leqslant C_6<\infty$.

Зазначимо, що для $\tilde{V}_{2^s}$ має місце рівність
$$%\begin{equation}\label{v-til}
\tilde{V}_{2^s}=\mu_{2^s}(\boldsymbol{x})=\prod\limits_{j=1}^d \mu_{2^s}(x_j),
$$%\end{equation}
де
$$
\mu_{2^s}(x_j)=\sqrt{\frac {\pi}{2}}
               \begin{cases}
                  1, & |x_j|\leqslant 2^s; \\
                  \frac{1}{2^s}(2^{s+1}-x_j), & 2^s<|x_j|\leqslant2^{s+1}; \\
                  0, & 2^{s+1}<|x_j|.
               \end{cases}
$$

Для функцій $g_1\in L_1(\mathbb{R}^d)$ та $g_2 \in  L_p(\mathbb{R}^d)$,  $1\leqslant p \leqslant \infty$, означимо їхню згортку згідно з формулою (див., наприклад, \cite[\S 1.5.1]{Nikolsky_1969_book})
$$
\big(g_1 \ast g_2\big)(\boldsymbol{x})=\frac{1}{(2\pi)^{\frac{d}{2}}}\int g_1(\boldsymbol{x}-\boldsymbol{u})g_2(\boldsymbol{u})d\boldsymbol{u}.
$$
При цьому виконується нерівність
\begin{equation}\label{zgortka_ner}
\|g_1\ast g_2\|_{L_p(\mathbb{R}^d)}\leqslant \frac{1}{(2\pi)^{\frac{d}{2}}} \|g_1\|_{L_1(\mathbb{R}^d)}\cdot\|g_2\|_{L_p(\mathbb{R}^d)}.
\end{equation}

Нехай $f\in L_p(\mathbb{R}^d)$, $1\leqslant p \leqslant \infty$. В такому випадку покладемо
\begin{equation}\label{sigma_2s}
\sigma_{2^s}(f,\boldsymbol{x})=\left(\frac{2}{\pi}\right)^{\frac{d}{2}}\big(V_{2^s} \ast f\big)(\boldsymbol{x})=\frac{1}{\pi^{d}}\int V_{2^s}(\boldsymbol{x}-\boldsymbol{u})f(\boldsymbol{u})d\boldsymbol{u}.
\end{equation}
Дана функція є аналогом періодичної суми Валле Пуссена порядку $2^s$, крім цього $\sigma_{2^s}(f,\boldsymbol{x})\in L_p(\mathbb{R}^d)$, $1\leqslant p \leqslant \infty$,  є цілою функцією експоненціального типу  $2^{s+1}$ по кожній змінній $x_j$, $j=\overline{1,d}$. У термінах перетворення Фур'є  $\sigma_{2^s}(f,\boldsymbol{x})$
можна подати у такому вигляді % \cite[ c.~359]{Nikolsky_1969_book}
%\begin{equation}\label{sigma_furie}
$$
\sigma_{2^s}(f,\boldsymbol{x})=\mathfrak{F}^{-1}(\mu_{2^s}\cdot\mathfrak{F} f).
$$
%\end{equation}

Далі, кожній функцій $f\in L_p(\mathbb{R}^d)$, $1\leqslant p \leqslant \infty$, поставимо у відповідність ряд
\begin{equation}\label{f_ryad}
f(\boldsymbol{x})=\sigma_{2^0}(f,\cdot)+\sum\limits_{s=1}^{\infty}\big(\sigma_{2^s}(f,\cdot)-\sigma_{2^{s-1}}(f,\cdot)\big),
\end{equation}
який збігається до функції $f$ у метриці простору $L_p(\mathbb{R}^d)$~\cite{Nikolskii_1968_sib}. Даний ряд будемо називати розкладом функції $f$ за сумами типу Валле Пуссена.

Введемо наступне позначення
\begin{equation}\label{qs}
q_0(f,\cdot)=\sigma_{2^0}(f,\cdot), \quad \quad q_{s}(f,\cdot)=\sigma_{2^s}(f,\cdot)-\sigma_{2^{s-1}}(f,\cdot), s\in \mathbb{N}.
\end{equation}
Тоді згідно з співвідношенням \eqref{qs} рівність \eqref{f_ryad} для $f$  можемо переписати у такому вигляді
$$
f(\boldsymbol{x})=\sum\limits_{s=0}^{\infty}q_s(f,\cdot),
$$
причому зауважимо, що носій перетворення Фур'є функції $q_s(f,\cdot)$ зосереджений на множині ${\Big\{\boldsymbol{\lambda}\colon 2^{s-1}\leqslant \max \limits_{j=\overline{1,d}}|\lambda_j|\leqslant 2^{s+1}\Big\}}$.

\begin{remark}\label{En-Vn}
Наближення функції $f\in L_p(\mathbb{R}^d)$, $1 \leqslant p \leqslant \infty$, за допомогою $\sigma_{2^s}(f,\boldsymbol{x})$ має такий же порядок, як і найкраще наближення цієї функції за допомогою функцій експоненціального типу $2^s$, тобто $E_{n}\big(f\big)_{L_p(\mathbb{R}^d)}$.
\end{remark}

Таким чином на основі твердження~\ref{T-decompoz}, відповідно до наведених позначень,  можна дати наступне означення просторів $B^{\Omega}_{p,\theta}(\mathbb{R}^d)$ (див. наслідок~1~\cite{Muronyuk_2013_zb}).

Функція $f$ належить простору $B^{\Omega}_{p,\theta}(\mathbb{R}^d)$, $1\leqslant p \leqslant \infty$ і $\Omega \in
\Phi_{\alpha, l}$, $\alpha>0$, $l\in \mathbb{N}$, якщо для неї скінченні величини
$$
\left(\sum
\limits_{s=0}^{\infty} \big(\Omega(2^{-s}))^{-\theta} \|q_{s}(f,\cdot)\|_{L_p(\mathbb{R}^d)}^{\theta}\right)^{\frac{1}{\theta}}, \ \mbox{при} \ \ 1\leqslant \theta < \infty
$$
та
$$
\sup \limits_{s\geqslant
0} \big(\Omega(2^{-s}))^{-1}\|q_{s}(f,\cdot)\|_{L_p(\mathbb{R}^d)}, \ \mbox{при}  \ \ \theta=\infty.
$$
При цьому норма $\|f\|_{B^{\Omega}_{p,\theta}(\mathbb{R}^d)}$, $1\leqslant \theta \leqslant \infty$, функцій $f$, згідно з теоремою~А,  задовольняє співвідношенням
\begin{equation}\label{Norm_dek}
\|f\|_{B^{\Omega}_{p,\theta}(\mathbb{R}^d)}\asymp\left(\sum
\limits_{s=0}^{\infty} \big(\Omega(2^{-s}))^{-\theta}  \|q_{s}(f,\cdot)\|^{\theta}_{L_p(\mathbb{R}^d)}\right)^{\frac{1}{\theta}},
\end{equation}
якщо $1\leqslant \theta < \infty$, i
\begin{equation}\label{Norm_infty_dek}
\|f\|_{B^{\Omega}_{p,\infty}(\mathbb{R}^d)}\asymp\sup \limits_{s\geqslant
0} \big(\Omega(2^{-s}))^{-1} \|q_{s}(f,\cdot)\|_{L_p(\mathbb{R}^d)}.
\end{equation}

Зазначимо, що як наслідок з твердження~\ref{T-decompoz} випливає також важливий для встановлення результатів
факт, що простори $B^{\Omega}_{p,\theta}(\mathbb{R}^d)$ зі зростанням параметра $\theta$
розширюються
\begin{equation}\label{vklad-Omega}
 B^{\Omega}_{p,1}(\mathbb{R}^d)\subset B^{\Omega}_{p,\theta}(\mathbb{R}^d) \subset B^{\Omega}_{p,\theta'}(\mathbb{R}^d) \subset
 B^{\Omega}_{p,\infty}(\mathbb{R}^d)\equiv H^{\Omega}_{p}(\mathbb{R}^d),  1\leqslant\theta<\theta' \leqslant
 \infty.
\end{equation}

Далі, якщо не стверджується інше, під терміном ``класи $B^{\Omega}_{p,\theta}(\mathbb{R}^d)$'' будемо розуміти
одиничні кулі у просторі $B^{\Omega}_{p,\theta}(\mathbb{R}^d)$, а саме
$$
B^{\Omega}_{p,\theta}(\mathbb{R}^d):= \left\{f\in L_p(\mathbb{R}^d) \colon \|f\|_{B^{\Omega}_{p,\theta}(\mathbb{R}^d)}\leqslant 1\right\}
$$
і при цьому збережемо для класів $B^{\Omega}_{p,\theta}B(\mathbb{R}^d)$ ті ж позначення, що і для
просторів $B^{\Omega}_{p,\theta}B(\mathbb{R}^d)$.

Тепер дамо означення апроксимаційної характеристики, яка буде досліджуватися у роботі.

Для $f\in L_q(\mathbb{R}^d)$, $1\leqslant q \leqslant \infty$, розглянемо частинну суму типу Валле Пуссена
\begin{equation}\label{Enf-Valle-Omega}
\mathbb{V}_n(f,\cdot)=\sum\limits_{s=0}^{n}q_s(f,\cdot)
\end{equation}
і покладемо
$$
\mathcal{E}_n\big(f\big)_{L_q(\mathbb{R}^d)}:=\big\|f(\cdot)-\mathbb{V}_{n-1}(f,\cdot)\big\|_{L_q(\mathbb{R}^d)}.
$$
Величина $\mathcal{E}_n\big(f\big)_{L_q(\mathbb{R}^d)}$ називається величиною наближенням функції $f$ частинними сумами типу Валле Пуссена. Якщо $F\subset L_q(\mathbb{R}^d)$, то  покладемо
\begin{equation}\label{En-F-Valle-Omega}
\mathcal{E}_n\big(F\big)_{L_q(\mathbb{R}^d)}:=\sup\limits_{f\in F} \mathcal{E}_n\big(f\big)_{L_q(\mathbb{R}^d)}.
\end{equation}

\section{Наближення частинними сумами типу Валле Пуссена у рівномірній та інтегральній метриках}

Попередньо сформулюємо твердження, яке буде істотно використовуватися при встановленні результатів.

\begin{proposition}[\!\!\cite{Nikolsky_1969_book}, \S\,3.3.4]\label{T-neriv}
Якщо
$1\leqslant p_1 \leqslant p_2 \leqslant \infty$, то для цілої функції
експоненціального типу $g=g_{\nu}\in L_p(\mathbb{R}^d)$ має місце
нерівність
\begin{equation}\label{Riz-Metric-Rd}
 \|g_{\nu}\|_{L_{p_2}(\mathbb{R}^d)}\leqslant 2^d\left( \prod \limits_{j=1}^d
 \nu_k\right)^{\frac{1}{p_1}-\frac{1}{p_2}}\|g_{\nu}\|_{L_{p_1}(\mathbb{R}^d)}.
\end{equation}
\end{proposition}

Нерівність \eqref{Riz-Metric-Rd} прийнято називати ``нерівністю різних метрик'' Нікольського для цілих функцій експоненціального типу.

Справедливе таке твердження.

\begin{theorem} \label{T-En-Om-iz-Rd}
Нехай $1\leqslant p<\infty$, $\Omega\in \Phi_{\alpha,l}$, де $\alpha> \frac{d}{p}$, $l\in \mathbb{N}$. Тоді для $1\leqslant \theta
\leqslant \infty$ має місце
порядкове співвідношення
\begin{equation} \label{T1-En-izotGer-infty}
   \mathcal{E}_n\big(B^{\Omega}_{p,\theta}(\mathbb{R}^d)\big)_{L_{\infty}(\mathbb{R}^d)}
   \asymp \Omega(2^{-n}) 2^{\frac{nd}{p}}.
\end{equation}
\end{theorem}

\begin{proof} Встановимо спочатку оцінку зверху в~\eqref{T1-En-izotGer-infty}. Оскільки, згідно з теоремою~А, для ${1\leqslant p<\infty}$ має місце вкладення $B^{\Omega}_{p,\theta}(\mathbb{R}^d)\subset H^{\Omega}_p(\mathbb{R}^d)$, $1\leqslant \theta < \infty$, то шукану оцінку достатньо отримати для величини $\mathcal{E}_{n}\big(H^{\Omega}_p(\mathbb{R}^d)\big)_{L_{\infty}(\mathbb{R}^d)}$. Далі для $f\in H^{\Omega}_p(\mathbb{R}^d)$, згідно з \eqref{Norm_infty_dek}, можемо записати $\|q_s(f,\cdot)\|_{L_p(\mathbb{R}^d)}\ll \Omega(2^{-s})$. Тому, використовуючи  нерівність Мінковського, нерівність різних метрик \eqref{Riz-Metric-Rd} та враховуючи, що $\Omega$ задовольняє мову $(S^{\alpha})$ з $\alpha> \frac{d}{p}$, отримуємо
$$
\mathcal{E}_n\big(B^{\Omega}_{p,\theta}(\mathbb{R}^d)\big)_{L_{\infty}(\mathbb{R}^d)}=
\Bigg\|f(\cdot)-\sum\limits_{s=0}^{n-1}q_s(f,\cdot)\Bigg\|_{L_{\infty}(\mathbb{R}^d)}=
\Bigg\|\sum\limits_{s=n}^{\infty}q_s(f,\cdot)\Bigg\|_{L_{\infty}(\mathbb{R}^d)}\ll
$$
$$
\ll \sum\limits_{s=n}^{\infty}2^{\frac{sd}{p}}\|q_s(f,\cdot)\|_{L_p(\mathbb{R}^d)}\ll \sum\limits_{s=n}^{\infty}2^{\frac{sd}{p}} \Omega(2^{-s})=
$$
$$
=\sum\limits_{s=n}^{\infty} \frac{\Omega(2^{-s})}{2^{-\alpha s}} 2^{-\left(\alpha-\frac{d}{p}\right)s} \ll
\frac{\Omega(2^{-n})}{2^{-\alpha n}} \sum\limits_{s=n}^{\infty} 2^{-\left(\alpha-\frac{d}{p}\right)s} \asymp
$$
$$
\asymp \frac{\Omega(2^{-n})}{2^{-\alpha n}} 2^{-\left(\alpha-\frac{d}{p}\right)n}= \Omega(2^{-n}) 2^{\frac{nd}{p}}.
$$

Оцінку зверху в~\eqref{T1-En-izotGer-infty} встановлено.

Перейдемо до встановлення  оцінки знизу в~\eqref{T1-En-izotGer-infty}. Оскільки згідно з \eqref{vklad-Omega} $B^{\Omega}_{p,1}(\mathbb{R}^d)\subset B^{\Omega}_{p,\theta}(\mathbb{R}^d)$, $1<\theta\leqslant \infty$, то шукану оцінку достатньо отримати для $\mathcal{E}_{n}\big(B^{\Omega}_{p,1}(\mathbb{R}^d)\big)_{L_{\infty}(\mathbb{R}^d)}$.

Нехай
$$
f_{n+1}(\boldsymbol{x})=\prod\limits_{j=1}^{d}\big(V_{2^{n+1}}(x_j)-V_{2^{n}}(x_j)\big), \ \ n\in\mathbb{N}, \ x_j\in \mathbb{R}, \ \ j=\overline{1,d}.
$$

У роботі~\cite{Yanchenko-UMG-2015} показано, що виконується оцінка
\begin{equation}\label{fn+1_norm_inf}
\|f_{n+1}(\cdot)\|_{L_{\infty}(\mathbb{R}^d)}\asymp 2^{nd},
\end{equation}
а у випадку $1\leqslant p< \infty$ в~\cite{Muronyuk_2013_zb} встановлено, що
\begin{equation}\label{fn+1_norm_all}
\|f_{n+1}(\cdot)\|_{L_p(\mathbb{R}^d)}\asymp 2^{nd\left(1-\frac{1}{p}\right)}, \ \ \ 1\leqslant p < \infty.
\end{equation}

Розглянемо тепер функцію
\begin{equation}\label{F1-izoGen}
F_1(\boldsymbol{x})=C_7 \Omega(2^{-n})2^{-nd\left(1-\frac{1}{p}\right)}f_{n+1}(\boldsymbol{x})
\end{equation}
і переконаємося, що при певному виборі сталої $C_7>0$, вона належить класу $B^{\Omega}_{p,1}(\mathbb{R}^d)$.

Оскільки носій перетворення Фур'є функції $f_{n+1}$ міститься на множині
$$
\Big\{\boldsymbol{\lambda}\colon 2^{n}\leqslant \max \limits_{j=\overline{1,d}}|\lambda_j|\leqslant 2^{n+2}\Big\},
$$
то, згідно з зазначеним вище, щодо носія перетворення Фур'є функцій $q_s(f,\cdot)$, ${s=0,1,\ldots}$,  отримаємо, що всі функції $q_s(f_{n+1},\cdot)$ окрім, можливо, $q_{n}(f_{n+1},\cdot)$, $q_{n+1}(f_{n+1},\cdot)$ та $q_{n+2}(f_{n+1},\cdot)$  тотожно дорівнюють нулеві.

Отже, згідно з  \eqref{Norm_dek} при $\theta=1$ будемо мати
$$
\|F_1(\cdot)\|_{B^{\Omega}_{p,1}(\mathbb{R}^d)}\asymp \sum\limits_{s\in\mathbb{Z}_+} \big(\Omega(2^{-s})\big)^{-1}\|q_s(F_1,\cdot)\|_{L_p(\mathbb{R}^d)}=
$$
$$
=C_7 \Omega(2^{-n})2^{-nd\left(1-\frac{1}{p}\right)} \sum\limits_{s\in\mathbb{Z}_+}\big(\Omega(2^{-s})\big)^{-1}\|q_s(f_{n+1},\cdot)\|_{L_p(\mathbb{R}^d)}=
$$
$$
=C_7 \Omega(2^{-n})2^{-nd\left(1-\frac{1}{p}\right)} \left(\frac{\|q_{n}(f_{n+1},\cdot)\|_{L_p(\mathbb{R}^d)}}{\Omega(2^{-n})} + \frac{\|q_{n+1}(f_{n+1},\cdot)\|_{L_p(\mathbb{R}^d)}}{\Omega(2^{-(n+1)})}+
\right.
$$
\begin{equation}\label{sum_qs}
\left. +\frac{\|q_{n+2}(f_{n+1},\cdot)\|_{L_p(\mathbb{R}^d)}}{\Omega(2^{-(n+2)})}\right).
\end{equation}

Для продовження оцінки \eqref{sum_qs} необхідно оцінити кожен з доданків $\|q_{s}(f_{n+1},\cdot)\|_{L_p(\mathbb{R}^d)}$ для $s=\{n, n+1, n+2\}$. Зауважимо, що враховуючи означення для $q_{s}(f_{n+1},\cdot)$ з урахуванням властивості згортки \eqref{zgortka_ner} і те, що згідно з властивістю 4) для $V_{2^s}$ норма ядра Валле Пуссена в $L_1(\mathbb{R}^d)$ є обмеженою, отримаємо, що  значення $\|q_{s}(f_{n+1},\cdot)\|_{L_p(\mathbb{R}^d)}$ буде залежати лише від значення норми $f_{n+1}$  в  $L_p(\mathbb{R}^d)$, а саме мають місце перетворення (див. також~\cite{Yanchenko-UMG-2015}):
$$
\|q_n(f_{n+1},\cdot)\|_{L_p(\mathbb{R}^d)}=\|\sigma_{2^n}(f_{n+1})-\sigma_{2^{n-1}}(f_{n+1})\|_{L_p(\mathbb{R}^d)} =
$$
$$
= \left(\frac{2}{\pi}\right)^{\frac{d}{2}}
\big\|\big(V_{2^n}(\cdot)-V_{2^{n-1}}(\cdot)\big)\ast f_{n+1}(\cdot)\big\|_{L_p(\mathbb{R}^d)}\leqslant
$$
$$
\leqslant\left(\frac{2}{\pi}\right)^{\frac{d}{2}} \frac{1}{(2\pi)^{\frac{d}{2}}}\|V_{2^n}(\cdot)-V_{2^{n-1}}(\cdot)\|_{L_1(\mathbb{R}^d)} \|f_{n+1}(\cdot)\|_{L_p(\mathbb{R}^d)}\leqslant C_8 2^{nd\left(1-\frac{1}{p}\right)}.
$$
Аналогічно отримуються оцінки і для $\|q_{n+1}(f_{n+1},\cdot)\|_{L_p(\mathbb{R}^d)}$ та $\|q_{n+2}(f_{n+1},\cdot)\|_{L_p(\mathbb{R}^d)}$.

Тоді \eqref{sum_qs} можемо продовжити таким чином
$$
\|F_1(\cdot)\|_{B^{\Omega}_{p,1}(\mathbb{R}^d)}\leqslant
$$
$$
\leqslant C_7 \Omega(2^{-n})2^{-nd\left(1-\frac{1}{p}\right)}  \left( C_8 \frac{2^{nd\left(1-\frac{1}{p}\right)}}{\Omega(2^{-n})}  + C_9  \frac{2^{nd\left(1-\frac{1}{p}\right)}}{\Omega(2^{-(n+1)})} +  C_{10} \frac{2^{nd\left(1-\frac{1}{p}\right)}}{\Omega(2^{-(n+2)})} \right) \leqslant C_{11},
$$
де $C_8, C_9, C_{10}, C_{11}>0$.

Отже, при відповідному виборі сталої $C_7$ функція $F_1$ належить класу $B^{\Omega}_{p,1}(\mathbb{R}^d)$.

Далі,  згідно з вибором функції $F_1$  і властивостями функції $f_{n+1}$ маємо $\mathbb{V}_{n-1}(F_1,\cdot)=0$. Тому, враховуючи \eqref{fn+1_norm_inf}, будемо мати
$$
\mathcal{E}_n\big(B^{\Omega}_{p,1}(\mathbb{R}^d)\big)_{L_{\infty}(\mathbb{R}^d)}\geqslant \mathcal{E}_n\big(F_1\big)_{L_{\infty}(\mathbb{R}^d)} =\|F_1(\cdot)-\mathbb{V}_{n-1}(F_1,\cdot)\|_{L_{\infty}(\mathbb{R}^d)}=
$$
$$
=\|F_1(\cdot)\|_{L_{\infty}(\mathbb{R}^d)}\asymp \Omega(2^{-n})2^{-nd\left(1-\frac{1}{p}\right)}\|f_{n+1}(\cdot)\|_{L_{\infty}(\mathbb{R}^d)}\asymp
$$
$$
\asymp \Omega(2^{-n})2^{-nd\left(1-\frac{1}{p}\right)} 2^{nd} \asymp \Omega(2^{-n}) 2^{\frac{dn}{p}}.
$$

Оцінки знизу отримано. Теорему~\ref{T-En-Om-iz-Rd} доведено.
\end{proof}

Зауважимо, що отримати подібне до теореми~\ref{T-En-Om-iz-Rd} твердження  можна і у випадку $1<p<\infty$. Для цього на основі твердження~\ref{T-decompoz} розглянемо декомпозиційне представлення для норми функцій з класів $B^{\Omega}_{p,\theta}(\mathbb{R}^d)$ в дещо іншій формі (див. наслідок~\cite{Muronyuk_2014_UMG}) і, відповідно, для оцінки знизу будемо розглядати певним чином видозмінену функцію $F_1$.

Спочатку розглянемо ядро Діріхле вигляду
$$
D_m (\boldsymbol{x})\prod \limits_{j=1}^{d}\frac{\sin m x_j} {x_j}, \quad m\in \mathbb{N}.
$$
Тоді для функції $f\in L_p(\mathbb{R}^d)$, $1<p<\infty$, покладемо
$$
S_{2^s}(f,\boldsymbol{x})=\frac{1}{\pi^d}\int
\limits_{\mathbb{R}^d} f(\boldsymbol{t}) D_{2^s}(\boldsymbol{x}-\boldsymbol{t}) d\boldsymbol{t},
$$
$$
f_{2^0}(\boldsymbol{x})=f_{(0)}(\boldsymbol{x})=S_{2^0}(f,\boldsymbol{x}),
$$
$$
f_{2^s}(\boldsymbol{x})=f_{(s)}(\boldsymbol{x})=S_{2^s}(f,\boldsymbol{x})-S_{2^{s-1}}(f,\boldsymbol{x}),  \quad \mbox{якщо} \quad
s \in \mathbb{N}.
$$
Крім того, у сенсі збіжності у метриці простору $L_p(\mathbb{R}^d)$, $1<p<\infty$, справджується рівність
\begin{equation}\label{f-Fourier-roz}
f(\boldsymbol{x}) = \sum\limits_{s=0}^{\infty} f_{(s)}(\boldsymbol{x}).
\end{equation}

Тоді на основі твердження~\ref{T-decompoz} функція $f(\boldsymbol{x})$ належить простору $B^{\Omega}_{p,\theta}(\mathbb{R}^d)$,
$r>0$, $1<p<\infty$, $1\leqslant \theta \leqslant \infty$, тоді і
тільки тоді, коли
$$
\|f\|_{B^{\Omega}_{p,\theta}(\mathbb{R}^d)}\asymp \left(\sum
\limits_{s=0}^{\infty}(\Omega(2^{-s}))^{-\theta}\|f_{(s)}(\cdot)\|^{\theta}_{L_p(\mathbb{R}^d)}\right)^{\frac{1}{\theta}}
< \infty,
$$
якщо $1\leqslant \theta < \infty$, i
\begin{equation}\label{Norm_infty_dek-fs}
\|f\|_{B^{\Omega}_{p,\theta}(\mathbb{R}^d)}\asymp \sup \limits_{s\geqslant
0} (\Omega(2^{-s}))^{-1}\|f_{(s)}(\cdot)\|_{L_p(\mathbb{R}^d)} < \infty,
\end{equation}
якщо $\theta=\infty$.

Далі для $f\in L_q(\mathbb{R}^d)$, $1<q<\infty$, розглянемо частинну суму \eqref{f-Fourier-roz} вигляду
\begin{equation}\label{Enf-Fourier-Omega}
S_n(f,\boldsymbol{x})=\sum \limits_{s=0}^{n}f_{(s)}(\boldsymbol{x})
\end{equation}
і, відповідно, наближення функції сумами~\eqref{Enf-Fourier-Omega}
\begin{equation}\label{En'-Fourier-Omega}
\mathcal{E}'_n\big(f\big)_{L_q(\mathbb{R}^d)}=\|f(\cdot)-S_n(f,\cdot)\|_{L_q(\mathbb{R}^d)}.
\end{equation}
Якщо ж $F\subset L_q(\mathbb{R}^d)$~--- деякий функціональний клас, то
\begin{equation}\label{En'-F-Fourier-Omega}
\mathcal{E}'_n\big(F\big)_{L_q(\mathbb{R}^d)}:=\sup\limits_{f\in F}\mathcal{E}'_n\big(f\big)_{L_q(\mathbb{R}^d)}.
\end{equation}

\begin{theorem} \label{T2-En-Om-iz-Rd}
Нехай $1< p<\infty$, $\Omega\in \Phi_{\alpha,l}$, де $\alpha> \frac{d}{p}$, $l\in \mathbb{N}$. Тоді для $1\leqslant \theta
\leqslant \infty$ має місце
порядкове співвідношення
\begin{equation} \label{T2-En-izotGer-infty}
   \mathcal{E}'_n\big(B^{\Omega}_{p,\theta}(\mathbb{R}^d)\big)_{L_{\infty}(\mathbb{R}^d)}
   \asymp \Omega(2^{-n}) 2^{\frac{nd}{p}}.
\end{equation}
\end{theorem}

\begin{proof} Оцінка зверху встановлюється цілком аналогічно, як і в теоремі~\ref{T-En-Om-iz-Rd}, але уже з використанням співвідношення~\eqref{Norm_infty_dek-fs}.

Для оцінки знизу, як уже відмічалося, побудуємо екстремальну функцію, яка її реалізує.  Нехай $\boldsymbol{k}=(k_1,\ldots,k_d) \in \mathbb{Z}^d_+$, розглянемо функцію~\cite{WangHeping-SunYongsheng-1995}
 $$
 D_{\boldsymbol{k}}(\boldsymbol{x})=\prod \limits_{j=1}^{d}D_{k_j}(x_j),
 $$
 де для кожного $k_j$ функція $D_{k_j}(x_j)$
 визначається таким чином:
 $$
 D_{k_j}(x_j)=\sqrt{\frac {2}{\pi}}\ \Big(2\sin {\frac
 {x_j}{2}}\cos{\frac {2k_j+1}{2}x_j} \Big) \cdot x_j^{-1}.
 $$

Далі покладемо
$$
 F_n(\boldsymbol{x})=\sum \limits_{k_1=2^n}^{2^{n+1}-1}\ldots \sum\limits_{k_d=2^n}^{2^{n+1}-1} D_{\boldsymbol{k}}(\boldsymbol{x})
$$
і зауважимо,  що для  $1<q\leqslant \infty$ справджується порядкова оцінка (див., наприклад, \cite{Muronyuk_2014_UMG, Yanchenko-Zb-2013})
\begin{equation}\label{D_k_2^n}
 \|F_n(\cdot)\|_{L_q(\mathbb{R}^d)}\asymp 2^{nd\left(1-\frac{1}{q}\right)}.
\end{equation}

Розглянувши функцію
$$
 F_2(\boldsymbol{x})=C_{12}  \Omega(2^{-n}) 2^{-nd\left(1-\frac{1}{p}\right)}F_n(\boldsymbol{x}), \quad C_{12}>0,
$$
 цілком аналогічно, як і в теоремі~\ref{T-En-Om-iz-Rd},  отримаємо оцінку \eqref{T2-En-izotGer-infty} у випадку, коли ${1<p<\infty}$.
\end{proof}

\begin{corollary} \label{Col-En-Om-iz-Rd}
Нехай $1< p<\infty$,  $r > \frac{d}{p}$. Тоді для $1\leqslant \theta
\leqslant \infty$ має місце
порядкове співвідношення
\begin{equation} \label{Col-En-izot-infty}
   \mathcal{E}'_n\big(B^{r}_{p,\theta}(\mathbb{R}^d)\big)_{L_{\infty}(\mathbb{R}^d)}
   \asymp  2^{-n\left(r-\frac{d}{p}\right)}.
\end{equation}
\end{corollary}

\begin{remark}
Відповідно до зауваження~\ref{En-Vn} точні за порядком оцінки величини \eqref{En-F-Valle-Omega}  для класів  $B^{\Omega}_{1,\theta}(\mathbb{R}^d)$ у випадку $1<q<\infty$ встановлено в~\cite{Muronyuk_2013_zb}, а для класів $B^{\Omega}_{p,\theta}(\mathbb{R}^d)$ у випадку $1<p\leqslant q<\infty$~--- в~\cite{Muronyuk_2014_UMG}. Також у згаданих роботах при цих же значеннях параметрів знайдено і оцінки величини $\mathcal{E}'_n\big(B^{\Omega}_{p,\theta}(\mathbb{R}^d)\big)_{L_q(\mathbb{R}^d)}$, крім того в~\cite{Muronyuk_2013_zb} знайдено оцінки і $\mathcal{E}'_n\big(B^r_{p,\theta}(\mathbb{R}^d)\big)_{L_q(\mathbb{R}^d)}$.
\end{remark}

Далі розглянемо випадок, коли параметри $p$ i $q$ в задачі про відшукання оцінки величини $\mathcal{E}_n\big(B^{\Omega}_{p,\theta}(\mathbb{R}^d)\big)_{L_q(\mathbb{R}^d)}$ приймають крайні значення, тобто $1$ або $\infty$.

\begin{theorem} \label{T3-En-Om-iz-Rd}
Нехай  $\Omega\in \Phi_{\alpha,l}$, де $\alpha>0$, $1\leqslant \theta
\leqslant \infty$. Тоді, для ${(p,q)=\big\{(1,1), (\infty,\infty)\big\}}$ справедливе
порядкове співвідношення
\begin{equation} \label{T2-En-izotGer-1infty}
   {\mathcal{E}}_n\big(B^{\Omega}_{p,\theta}(\mathbb{R}^d)\big)_{L_q(\mathbb{R}^d)}
   \asymp \Omega(2^{-n}).
\end{equation}
\end{theorem}

\begin{proof} Оскільки має місце вкладення \eqref{vklad-Omega}, то оцінку зверху достатньо встановити для класів $B^{\Omega}_{p,\infty}(\mathbb{R}^d)\equiv H^{\Omega}_p(\mathbb{R}^d)$. Як уже зазначалося, для $f\in H^{\Omega}_p(\mathbb{R}^d)$ з \eqref{Norm_infty_dek} випливає співвідношення $\|q_s(\cdot)\|_{L_q(\mathbb{R}^d)}\ll  \Omega(2^{-s})$. Тоді, скориставшись нерівністю Мінковського, можемо записати
$$
{\mathcal{E}}_n\big(f\big)_{L_p(\mathbb{R}^d)}\leqslant\Bigg\|f(\cdot)-\sum\limits_{s=0}^{n-1}q_s(\cdot)\Bigg\|_{L_p(\mathbb{R}^d)}=
\Bigg\|\sum\limits_{s=n}^{\infty}q_s(\cdot)
\Bigg\|_{L_p(\mathbb{R}^d)}\ll
$$
\begin{equation} \label{Omega-1infty}
\ll \sum\limits_{s=n}^{\infty}\|q_s(\cdot)\|_{L_q(\mathbb{R}^d)}\leqslant \sum\limits_{s=n}^{\infty} \Omega(2^{-s}).
\end{equation}
Далі, врахувавши, що $\Omega$ задовольняє мову $(S^{\alpha})$ з $\alpha> 0$, оцінку \eqref{Omega-1infty} продовжимо таким чином
$$
\sum\limits_{s=n}^{\infty} \Omega(2^{-s})= \sum\limits_{s=n}^{\infty} \frac{\Omega(2^{-s})}{2^{-\alpha s}} 2^{-\alpha s} \ll
\frac{\Omega(2^{-n})}{2^{-\alpha n}} \sum\limits_{s=n}^{\infty} 2^{- \alpha s} \asymp
$$
$$
\asymp \frac{\Omega(2^{-n})}{2^{-\alpha n}} 2^{- \alpha n}= \Omega(2^{-n}).
$$

Оцінку зверху в \eqref{T2-En-izotGer-1infty} встановлено.

Оцінку знизу в \eqref{T2-En-izotGer-1infty}, як уже відмічалося при доведенні теореми~\ref{T-En-Om-iz-Rd}, достатньо встановити для класу $B^{\Omega}_{p,1}(\mathbb{R}^d)$. Для цього розглянемо функції:
$$
F_3(\boldsymbol{x})=C_{13} \Omega(2^{-n})2^{-nd}f_{n+1}(\boldsymbol{x}), \ \ C_{13}>0, \ \mbox{коли} \ p=\infty
$$
і
$$
F_4(\boldsymbol{x})=C_{14} \Omega(2^{-n}) f_{n+1}(\boldsymbol{x}), \ \ C_{14}>0, \ \mbox{коли} \ p=1.
$$

Покажемо, що  дані функції належать класу $B^{\Omega}_{p,1}(\mathbb{R}^d)$ при певному виборі сталих $C_{13}$ та $C_{14}$.

Для $F_3$, міркуючи аналогічно як і в теоремі~\ref{T-En-Om-iz-Rd} для $F_1$, скориставшись оцінкою \eqref{fn+1_norm_inf}, маємо:
$$
\|F_3(\cdot)\|_{B^{\Omega}_{\infty,1}(\mathbb{R}^d)}\asymp \sum\limits_{s\in\mathbb{Z}_+} \big(\Omega(2^{-s})\big)^{-1}\|q_s(F_3,\cdot)\|_{L_{\infty}(\mathbb{R}^d)}=
$$
$$
=C_{13} \Omega(2^{-n})2^{-nd} \sum\limits_{s\in\mathbb{Z}_+}\big(\Omega(2^{-s})\big)^{-1}\|q_s(f_{n+1},\cdot)\|_{L_{\infty}(\mathbb{R}^d)}=
$$
$$
=C_{13} \Omega(2^{-n})2^{-nd} \left(\frac{\|q_{n}(f_{n+1},\cdot)\|_{L_{\infty}(\mathbb{R}^d)}}{\Omega(2^{-n})} + \frac{\|q_{n+1}(f_{n+1},\cdot)\|_{L_{\infty}(\mathbb{R}^d)}}{\Omega(2^{-(n+1)})}+
\frac{\|q_{n+2}(f_{n+1},\cdot)\|_{L_{\infty}(\mathbb{R}^d)}}{\Omega(2^{-(n+2)})}\right)\leqslant
$$
$$
\leqslant C_{13} \Omega(2^{-n})2^{-nd} \left(C_{15} \frac{2^{nd}}{\Omega(2^{-n})} + C_{16}  \frac{2^{nd}}{\Omega(2^{-(n+1)})} +  C_{17} \frac{2^{nd}}{\Omega(2^{-(n+2)})} \right) \leqslant C_{18},
$$
де $C_{15}, C_{16}, C_{17}, C_{18} >0$.

Відповідно для $F_4$, з урахуванням того, що при $p=1$, згідно з \eqref{fn+1_norm_all} маємо $\|f_{n+1}(\cdot)\|_{L_1(\mathbb{R}^d)}\leqslant C_{19}$, $C_{19}>0$, отримаємо
$$
\|F_4(\cdot)\|_{B^{\Omega}_{1,1}(\mathbb{R}^d)}\asymp \sum\limits_{s\in\mathbb{Z}_+} \big(\Omega(2^{-s})\big)^{-1}\|q_s(f_{n+1},\cdot)\|_{L_1(\mathbb{R}^d)}\leqslant C_{20}, \quad C_{20}>0.
$$

Окрім того, згідно з побудовою,  $\mathbb{V}_{n-1}(F_3,\cdot)=0$ і $\mathbb{V}_{n-1}(F_4,\cdot)=0$. Тому, враховуючи \eqref{fn+1_norm_inf}, будемо мати:
$$
\mathcal{E}_n\big(B^{\Omega}_{\infty,1}(\mathbb{R}^d)\big)_{L_{\infty}(\mathbb{R}^d)}\geqslant \mathcal{E}_n\big(F_3\big)_{L_{\infty}(\mathbb{R}^d)} =\|F_3(\cdot)-\mathbb{V}_{n-1}(F_3,\cdot)\|_{L_{\infty}(\mathbb{R}^d)}=
$$
$$
=\|F_2(\cdot)\|_{L_{\infty}}\asymp \Omega(2^{-n})2^{-nd}\|f_{n+1}(\cdot)\|_{L_{\infty}}\asymp \Omega(2^{-n}).
$$
Аналогічно, згідно з \eqref{fn+1_norm_all}, для $F_4$ отримуємо
$$
\mathcal{E}_n\big(B^{\Omega}_{1,1}(\mathbb{R}^d)\big)_{L_1(\mathbb{R}^d)}\geqslant \mathcal{E}_n(F_4)_{L_1(\mathbb{R}^d)} =\|F_4(\cdot)\|_{L_1(\mathbb{R}^d)}\asymp \Omega(2^{-n}).
$$

Оцінки знизу встановлено. Теорему~\ref{T3-En-Om-iz-Rd} доведено.
\end{proof}

Зробимо декілька коментарів стосовно одержаних результатів.

Проаналізувавши одержані оцінки  \eqref{T1-En-izotGer-infty}, \eqref{T2-En-izotGer-infty} і  \eqref{T2-En-izotGer-1infty} бачимо, що вони не залежать від значення параметра $\theta$, тобто є однаковими як для узагальнених класів типу Нікольського, так і для узагальнених класів типу Бєсова.

Для $\Omega(t)=t^r$, $0<r<l$, тобто для класів $B^{r}_{p,\theta}(\mathbb{R}^d)$, відповідні до теорем~\ref{T-En-Om-iz-Rd} і \ref{T3-En-Om-iz-Rd} результати встановлено в~\cite{Yanchenko-UMG-2015}. Ізотропні класи Нікольського--Бєсова  $B^{r}_{p,\theta}(\mathbb{R}^d)$ також досліджувалися в~\cite{Yanchenko_CMP2021}.

Як уже зазначалося в одновимірному випадку ($d=1$) класи $B^{\Omega}_{p,\theta}(\mathbb{R})$ збігаються з класами  $S^{\Omega}_{p,\theta}B(\mathbb{R})$, відповідно оцінка теореми~\ref{T2-En-Om-iz-Rd} раніше встановлена в~\cite{Yanchenko-Zb-2014}, а наслідку~\ref{Col-En-Om-iz-Rd}~--- в \cite{Yanchenko-Zb-2013}. Результати теорем~\ref{T-En-Om-iz-Rd} і \ref{T3-En-Om-iz-Rd} є новими і в одновимірному випадку.

Відмітимо також, що ізотропні класи Нікольського--Бєсова  періодичних функцій багатьох змінних та їхні узагальнення з точки зору знаходження точних за порядком оцінок різних апроксимаційних характеристик досліджувалися, зокрема, у роботах \cite{Derevyanko-2014-UMG, Myronyuk-UMG-2012-9, Romanyuk-UMG-2009, Romanyuk-2008-zb, RomanyukAS-2013MN, Romanyuk-AV-UMG-2009, Romanyuk-other-2023CMP, Stasyuk-2011-MS}.

\vskip 5mm

\section{Acknowledgments}

\emph{This research is partially supported by the Volkswagen Foundation project ``From Modeling
and Analysis to Approximation'' and the Grant H2020-MSCA-RISE-2019, project number
873071 (SOMPATY: Spectral Optimization: From Mathematics to Physics and Advanced
Technology).}

%%%%%%%%%%%%%%%%%%%%%%%%%%%%%%%%%%%%%%%%%%%%%%%%%%

\vskip 25 mm

CONTACT INFORMATION

\bigskip

\begin{minipage}{7.5 cm}
M.\,I.~Gromyak \\
Ternopil Volodymyr Hnatiuk \\ National Pedagogical University, \\
2 Maxyma Kryvonosa str., Ternopil, \\
46027, Україна\\
ghromjak@tnpu.edu.ua \\
\end{minipage}
\begin{minipage}{7.5 cm}
М.\,І.~Гром'як\\
Тернопільський національний \\ педагогічний  університет \\
імені Володимира Гнатюка, \\
вул. Максима Кривоноса 2, Тернопіль, \\
46027, Україна\\
ghromjak@tnpu.edu.ua
\end{minipage}

\bigskip

\begin{minipage}{7.5 cm}
O.\,Ya.~Radchenko \\
Ternopil Volodymyr Hnatiuk \\ National Pedagogical University, \\
2 Maxyma Kryvonosa str., Ternopil, \\
46027, Україна\\
yan.olga1208@gmail.com \\
\end{minipage}
\begin{minipage}{7.5 cm}
O.\,Я.~Радченко\\
Тернопільський національний \\ педагогічний  університет \\
імені Володимира Гнатюка, \\
вул. Максима Кривоноса 2, Тернопіль, \\
46027, Україна\\
yan.olga1208@gmail.com
\end{minipage}

\bigskip

\begin{minipage}{7.5 cm}
S.\,Ya.~Yanchenko \\
Institute of Mathematics of NASU  \\
3 Tereschenkivska str.,  Kyiv, \\
01024, Україна\\
yan.sergiy@gmail.com
\end{minipage}
\begin{minipage}{7.5 cm}
С.\,Я.~Янченко\\
Інститут математики НАН України, \\
вул. Терещенківська 3, Київ, \\
01024, Україна\\
yan.sergiy@gmail.com
\end{minipage}

\end{document}